\documentclass[10pt]{amsart}

\usepackage{graphicx}
\usepackage{amssymb}	
\usepackage{amsmath}
\usepackage[ruled,vlined]{algorithm2e} 
\usepackage[hidelinks]{hyperref}

\title{
Optimal mode decomposition for control
} 
  
\author{Lucas Mieg}
\author{Martin M\"onnigmann}%
\thanks{\url{https://git.noc.rub.de/rub-control/omdc}}%
\thanks{This work was funded by the Deutsche Forschungsgemeinschaft (DFG, German Research Foundation) - Project-ID 422037413 - TRR 287.
Lucas Mieg and Martin M\"onnigmann are with Automatic Control and Systems Theory, Department of Mechanical Engineering, Ruhr-Universit\"at Bochum, Germany
{\tt\small \{lucas.mieg, martin.moennigmann\}@ruhr-uni-bochum.de}}%

\begin{document}

\begin{abstract}
We present an extension of \textit{optimal} mode decomposition (OMD) for autonomous systems to systems with controls. 
The extension is developed along the same lines as the extension of \textit{dynamic} mode decomposition (DMD) to DMD with control (DMDc).

DMD identifies a linear dynamic system from high-dimen\-sional snapshot data.
DMD is often combined with a subsequent reduction by a projection to a truncated basis for the space spanned by the snapshots. 
In OMD, the identification and reduction are essentially integrated into a single optimization step, thus avoiding the somewhat adhoc decoupled, a posteriori reduction that is necessary if DMD is to be used for model reduction.
DMD was devised for autonomous systems and later extended to DMD for systems with control inputs (DMDc).
We present the analogous extension of OMD to OMDc, i.e.\ OMD for systems with control inputs. We illustrate the proposed method with an application to coupled diffusion-equations that model the drying of a wood chip. 
Reduced models of this type are required for the efficient simulation of industrial drying processes. 
\end{abstract}

\maketitle

\section{Introduction}\label{sec:Introduction}

Starting with the introduction of dynamic mode decomposition (DMD) by~\cite{SCHMID_2010},
the analysis of linear characteristics in flow fields, or distributed systems in general, has gained a lot of attention. 
Because DMD is often combined with a projection for model reduction, it also has become an alternative to classical model reduction techniques. 
Here ``classical techniques'' refers to methods based on proper orthogonal decomposition (POD) and projection that were originally developed to analyse coherent structures in turbulent flows~\cite{sirovich1987coherent}.
DMD and related techniques are an attractive alternative because they involve operations on only the data that describes the process. 
In contrast, the classical techniques involve a Galerkin (or Petrov-Galerkin) projection of the Navier-Stokes equations that govern the flow field, or, more generally, the partial differential equations that govern the system at hand. 

We emphasize that DMD generates an approximation to the Koopman operator of the dynamical system.
DMD has been adapted to include nonlinear mappings of the original data, leading to the so-called \textit{extended} DMD (eDMD)~\cite{williams2015edmd}.
The present paper is restricted to  optimal linear projections of the data.

DMD delivers two components: a set of orthonormal modes is identified, and a linear model for the dynamics is found such that each mode is equipped with an eigenvalue, i.e. a frequency of oscillation and a rate of decay.
The resulting system is linear but typically large, as the system dimension is equal to the number of linearly independent snapshots. 
If DMD is used for model reduction, the set of modes is truncated. 
It was shown in~\cite{wynn2013optimal} that the truncation of DMD modes analogously to POD modes delivers suboptimal models. As a remedy, \textit{optimal} mode decomposition (OMD) incorporates a rank constraint on the truncated basis into the optimization problem that DMD is based on~\cite{goulart2012optimal}. 
In other words, OMD determines, given the desired number of basis vectors, the optimal truncated basis, while DMD first finds the optimal basis which is then truncated subsequently.
The optimization with the additional rank constrained can be cast as an optimization on the Grassmann manifold~\cite{goulart2012optimal}. 
A mixed-norm optimization has been proposed as an alternative~\cite{jovanovic2014sparsity}, where a 1-norm regularization term promotes a sparse set of modes, which, however, is not optimal in the OMD sense.
Both OMD and sparsity-promoting DMD have been developed for autonomous systems. 

We first formalize the system identification problem and review DMDc as needed for our paper in section~\ref{sec:dmdc}. 
The extension of OMD to OMDc is presented in section~\ref{sec:subspace-optimization}. 
Remarks on the computationally efficient implementation are also given in section~\ref{sec:subspace-optimization}. 
Section~\ref{sec:woodchip} presents results for an application of OMDc to coupled partial differential equations that describe the drying of wet particles, which are used in models of industrial wood drying processes. 
Conclusions are stated in section~\ref{sec:conclusions}. 

\section{Dynamic mode Decomposition for Control}
\label{sec:dmdc}
Dynamic mode decomposition for control (DMDc~\cite{proctor2016dynamic}) finds a linear approximation to the mapping 
\begin{equation*}
  \mathcal{F}:\mathcal{X}\times\mathcal{U}\rightarrow\mathcal{X}
   : (x_k,\,u_k) \rightarrow x_{k+1},
\end{equation*}
where $\mathcal{X}\subset\mathbb{R}^n$, $\mathcal{U}\subset\mathbb{R}^{p}$, $x_k$ denotes the system state at time $t_k=k\Delta t$, and $u_k$ denotes the system inputs at $t_k$.
We stress that the number of states $n$ is usually much larger than the number of time steps $m$.
This is, for example, the case for finite-volume simulation results, where the number of discrete volumes $n$ is usually much larger than the number of snapshots $m$.

DMDc approximates $\mathcal{F}$ with a time-discrete linear system
\begin{equation}
  \label{eq:discrete-lti}
  x_{k+1} = A x_k +  B u_k,
  \,\, A\in\mathbb{R}^{n\times n}, \, B \in\mathbb{R}^{n\times p}.
\end{equation}
In its first step, DMDc collects the system states $x_k$ 
in the columns of the snapshot matrix $S$
\begin{equation}
  S = [\cdots x_k \cdots] \in \mathbb{R}^{n\times m}, \quad k = 0,\ldots,m-1.
\end{equation}
After splitting the snapshot matrix $S$ into
\begin{equation}
  \label{eq:def-snaps}
  X = [x_0 \cdots x_{m-2}], \,\,
  Y = [x_1 \cdots x_{m-1}]
\end{equation}
the problem of finding the optimal matrices $A$ and $B$ for~\eqref{eq:discrete-lti} can be stated as 
\begin{equation}
  \label{eq:def-dmdc-prob}
  \min_{A,B} \Vert Y - AX - BU \Vert_\text{F},
\end{equation}
where $\Vert\cdot\Vert_\text{F}$ denotes the Frobenius norm, and $U\in\mathbb{R}^{p\times m-1}$ collects the sequence of inputs that yielded $S$ in its columns.
With the abbreviations
\begin{equation}\label{eq:HelperOmegaAndQ}
  \Omega = \begin{bmatrix}
    X \\ U
  \end{bmatrix}, \quad
  G = \begin{bmatrix}
    A \,\, B
  \end{bmatrix}
\end{equation}
problem~\eqref{eq:def-dmdc-prob} can be rewritten as a least-squares problem of the form
\begin{equation}\label{eq:StandardLeastSquares}
  \min_{G} \Vert Y - G \Omega \Vert_\text{F}
\end{equation}
which can be solved with the pseudo-inverse~\cite{proctor2016dynamic}
\begin{equation}\label{eq:PseudoInv}
  \Omega^\dagger = \tilde{V} \tilde{\Sigma}^{-1} \begin{bmatrix}
    \tilde{\Phi}_1^\intercal \,\,\, \tilde{\Phi}_2^\intercal,
  \end{bmatrix}  
\end{equation}
where $\tilde{V}\in\mathbb{R}^{(m-1)\times d}$, $\tilde{\Sigma}\in\mathbb{R}^{d\times d}$,
$\tilde{\Phi}_1\in\mathbb{R}^{n\times d}$ 
and $\tilde{\Phi}_2\in\mathbb{R}^{p\times d}$ with $d=\text{rank}\, \Omega \le m-1$
result from the thin singular value decomposition
$\begin{bmatrix}
      \tilde{\Phi}_1 \\ \tilde{\Phi}_2
    \end{bmatrix} \tilde{\Sigma} \tilde{V}^\intercal = \Omega$
~\cite[§2.4.3]{golub2013matrix}.
Note that the row dimensions of $\tilde{\Phi}_1$ and $\tilde{\Phi}_2$ correspond to those of $X$ and $U$, respectively.
Right-multiplying the objective function of~\eqref{eq:StandardLeastSquares} by $\Omega^\dagger$ from~\eqref{eq:PseudoInv} yields the optimal approximation $[\tilde{A}\,\tilde{B}]$ of $G= \left[A\, B\right]$ introduced in~\eqref{eq:HelperOmegaAndQ}
as $[\tilde{A}\, \tilde{B}]= Y\,\Omega^\dagger,$
or equivalently, 
\begin{equation}
  \label{eq:dmdc-matrices-full}
  \begin{aligned}
    \tilde{A} &= Y \tilde{V} \tilde{\Sigma}^{-1}\tilde{\Phi}_1^\intercal \\
    \tilde{B} &= Y \tilde{V} \tilde{\Sigma}^{-1}\tilde{\Phi}_2^\intercal.
  \end{aligned}
\end{equation}

We recall the dimension $n$ of the state space is usually large, because $n$ is the number of control volumes in a finite-volume simulation, for example. Consequently, $\tilde{A}\in\mathbb{R}^{n\times n}$ in~\eqref{eq:dmdc-matrices-full} may be prohibitively large. 
As a remedy, $\tilde{A}$ can be approximated with a second thin singular value decomposition. 
The decomposition $Y=\hat{\Phi}\hat{\Sigma}\hat{V}^\intercal$ 
is used for this purpose, which provides an approximation of $\tilde{A}$ and $\tilde{B}$ in the column space of $Y$, i.e.\ the output space~\cite{proctor2016dynamic}.
More precisely, 
let $\hat{\Phi}_r \hat{\Sigma}_r \hat{V}^\intercal_r$ denote the truncation of the singular value decomposition of $Y$ to its $r$ largest singular values with $r\le \text{rank}\, Y\le m-1$. 
Then
\begin{equation}
  \begin{aligned}
    \hat{A} &= \hat{\Phi}^\intercal_r \tilde{A} \hat{\Phi}_r
     &=& \hat{\Sigma}_r \hat{V}^\intercal_r \tilde{V} \tilde{\Sigma}^{-1}\tilde{\Phi}_1^\intercal  \hat{\Phi}_r
     \,&\in\mathbb{R}^{r\times r},
\\
    \hat{B} &= \hat{\Phi}^\intercal_r \tilde{B}
     &=& \hat{\Sigma}_r \hat{V}^\intercal_r \tilde{V} \tilde{\Sigma}^{-1}\tilde{\Phi}_2^\intercal
    \,&\in\mathbb{R}^{r \times p}\,\,
  \end{aligned}
\end{equation}
provide the approximations 
$\tilde{A} \approx \hat{\Phi}_r \hat{A} \hat{\Phi}^\intercal_r$ and 
$\tilde{B} \approx \hat{\Phi}_r \hat{B}$ with leading dimension $r< n$. 

Although the truncation $\hat{\Phi}_r \hat{\Sigma}_r \hat{V}^\intercal_r$ is optimal for approximating $Y$ in the sense of Eckhart-Young's Theorem~\cite[Th.~2.4.8]{golub2013matrix}, 
it is somewhat arbitrary to first find $A$ and $B$ by minimizing their approximation error with~\eqref{eq:def-dmdc-prob}, and to only subsequently determine a reduced state space, 
thus introducing a second approximation. 
In other words, the $r$-dimensional subspace spanned by the columns of $\hat{\Phi}_r$, is in general not the optimal $r$-dimensional subspace for solving~\eqref{eq:def-dmdc-prob}.
It is the central idea of the \textit{optimal} mode decomposition proposed in~\cite{goulart2012optimal} to replace the \textit{dynamic} mode decomposition with subsequent projection explained so far by an approach that accounts for both steps simultaneously.
In~\cite{goulart2012optimal}, this \textit{optimal} mode decomposition was introduced for autonomous systems. 
We extend the approach to systems with inputs in the next section. Technically this is done by accounting for a projection matrix $L\in\mathbb{R}^{n\times r}$ in~\eqref{eq:def-dmdc-prob}, instead of using an adhoc defined $\hat{\Phi}_r$. 

\section{Optimal mode decomposition for control}
\label{sec:subspace-optimization}
The goal of optimal mode decomposition for control is to find matrices 
$M\in\mathbb{R}^{r\times r}$, $P\in\mathbb{R}^{r\times p}$, and $L\in\mathbb{R}^{n\times r}$ such that the states of the system~\eqref{eq:discrete-lti} can be approximated as
\begin{subequations}\label{eq:RedSys}
  \begin{align}
    a_{k+1} &= M a_k + P u_k \\
    x_k &\approx L a_k
  \end{align}
\end{subequations}
where~(\ref{eq:RedSys}a) replaces the $n$-dimensional dynamical system~\eqref{eq:discrete-lti} with a lower-dimen\-sional system and state $a_k\in\mathbb{R}^r$, and $L$ lifts the low-dimensional state to the original state dimension in~(\ref{eq:RedSys}b). 
We follow~\cite{goulart2012optimal} in using $L^\intercal$ for projection and assuming $L^\intercal L= I$. Substituting the projection operation $a_k= L^\intercal x_k$ into~(\ref{eq:RedSys}a) and multiplying by $L$ results in
\begin{equation}\label{eq:FomHelper}
  x_{k+1}= LML^\intercal x_k+ LPu_k. 
\end{equation}
This yields $LML^\intercal= A$ when compared to~\eqref{eq:discrete-lti}, 
which illustrates that $L^\intercal$ acts as a projector from $\mathbb{R}^n$ to $\mathbb{R}^r$, $M$ evolves the $r$-dimensional system, and $L$ lifts back to $\mathbb{R}^n$. 

By replacing the original dynamics~\eqref{eq:discrete-lti} with \eqref{eq:FomHelper} in the \textit{dynamic} mode decomposition optimization problem~\eqref{eq:omdc-prob-full}, the \textit{optimal} mode decomposition optimization problem
\begin{subequations}
  \label{eq:omdc-prob-full}
  \begin{align}
    \label{eq:def-omdc-prob}
    \min_{L,M,P} \Vert & Y - L M L^\intercal X - L P U \Vert_\mathrm{F}^2 \\
    \label{eq:constraint-l}
    \quad &\mathrm{s.t.} \quad L^\intercal L = I
  \end{align}  
\end{subequations}
results.

We follow~\cite{goulart2012optimal} in solving~\eqref{eq:omdc-prob-full} in two steps.
The system and input matrix are first eliminated such that a problem in $L$ remains.
Subsequently, the constraint~\eqref{eq:constraint-l} is treated by optimizing on the Grassmann manifold. 

\subsection{Elimination of system matrices}
Assume the mode matrix $L$ to be arbitrary but fixed first.
In this case the system matrix $M$ and the input matrix $P$ can be eliminated with the first order necessary optimality condition for~\eqref{eq:def-omdc-prob}.
Let $F(L,M,P)$ refer to the cost function in~\eqref{eq:def-omdc-prob} and 
recall that the Frobenius norm of some matrix $C$ can be expressed with the trace
\begin{equation*}
  \Vert C \Vert_F^2 = \sum_{i,j} C_{ij}C_{ij} = \text{tr}\,C^\intercal C
\end{equation*}
For the cost function $F(L,M,P)$, this yields
\begin{equation*}
  \begin{aligned}
    F(&L,M,P) = \mathrm{tr}\Bigl( \bigl(Y - L M L^\intercal X\bigr)^\intercal \bigl(Y - L M L^\intercal X\bigr)\Bigr) \\
    -2 &\mathrm{tr}\Bigl(\bigl(Y - L M L^\intercal X\bigr)^\intercal\bigl(LPU\bigr)\Bigr) 
     + \mathrm{tr}\bigl(U^\intercal P^\intercal P U\bigr)
  \end{aligned}
\end{equation*}
Solving the necessary condition for optimality 
\begin{equation}
  \frac{\partial}{\partial P} F = 2 P U U^\intercal - 2 L^\intercal \big( Y - L M L^\intercal X \big) U^\intercal= 0
\end{equation}
yields
\begin{equation}
  \label{eq:optimal-p}
  P^\ast = L^\intercal Y U^\intercal \bigl(U U^\intercal\bigr)^{-1} - M L^\intercal X U^\intercal \bigl(U U^\intercal\bigr)^{-1} .
\end{equation}
After substituting $P=P^\ast$ into~\eqref{eq:def-omdc-prob}, it remains to find the minimizing $M$ and $L$, i.e.
\begin{equation}
  \begin{aligned}
    \min_{M,L} & \Vert \bigl(I - L L^\intercal\bigr) Y + L L^\intercal \hat{Y} - L M L^\intercal \hat{X} \Vert_\mathrm{F}^2 \\
    \hat{Y} = Y Q,& \quad \hat{X} = X Q, \quad
    Q = I - U^\intercal \bigl(U U^\intercal\bigr)^{-1}U,
  \end{aligned}
\end{equation}
where $Q\in\mathbb{R}^{(m-1)\times(m-1)}$ is an orthogonal projector by definition, which implies $Q Q = Q$ and $Q^\intercal=Q$.
Furthermore, note that $\bigl(I - L L^\intercal\bigr)L = 0$ for any $L$,
which simplifies 
the cost function $F$ to 
\begin{equation*}
  \begin{aligned}
    F(L,M,P^\ast) = & \, \mathrm{tr}\Bigl( Y^\intercal \bigl(I - L L^\intercal\bigr) Y \Bigr) \\
    -2 & \,\mathrm{tr} \Bigl(\hat{Y}^\intercal L M L ^\intercal \hat{X}\Bigr) 
     + \, \mathrm{tr} \Bigl(\hat{X}^\intercal L M^\intercal M L^\intercal \hat{X}\Bigr).
  \end{aligned}
\end{equation*}
The partial derivative of the cost function $F$ with respect to $M$ must fulfill the necessarily condition for optimality, i.e. $\frac{\partial}{\partial M}F=0$, which yields
\begin{equation}
    M^\ast = L^\intercal \hat{Y} \hat{X}^\intercal L \bigl(L^\intercal\hat{X}\hat{X}^\intercal L\bigr)^{-1}.
\end{equation}
Finally, the optimization problem that remains for $L$ reads
\begin{equation}
  \label{eq:omdc-cost-l}
  \begin{aligned}
    \min_{L} \Vert \bigl(I -& L L^\intercal\bigr) Y + L L^\intercal \hat{Y} \\
    &- L L^\intercal \hat{Y} \hat{X}^\intercal L \bigl(L^\intercal\hat{X}\hat{X}^\intercal L\bigr)^{-1} L^\intercal \hat{X} \Vert_\mathrm{F}^2    
  \end{aligned}
\end{equation}
In this case, the derivative required for the optimality conditions reads
\begin{equation}
  \label{eq:cwise-derivative}
  \begin{aligned}
    \frac{\partial}{\partial L} & F = -2 Y U^\intercal \bigl(UU^\intercal\bigr)^{-1} U Y^\intercal L \\
    -&2 \hat{Y}\hat{X}^\intercal L \bigl(L^\intercal\hat{X}\hat{X}^\intercal L\bigr)^{-1} L^\intercal \hat{X} \hat{Y}^\intercal L \\
    -&2 \hat{X}\hat{Y}^\intercal L  L^\intercal \hat{Y} \hat{X}^\intercal L \bigl(L^\intercal\hat{X}\hat{X}^\intercal L\bigr)^{-1} \\
    +&2 \hat{X}\hat{X}^\intercal L \bigl(L^\intercal\hat{X}\hat{X}^\intercal L\bigr)^{-1} L^\intercal \hat{Y} \hat{X}^\intercal L \\
    & \qquad \qquad \qquad \qquad \quad L^\intercal \hat{X} \hat{Y}^\intercal L \bigl(L^\intercal\hat{X}\hat{X}^\intercal L\bigr)^{-1}.
  \end{aligned}
\end{equation}

In contrast to the previous steps, the necessary condition of optimality~\eqref{eq:cwise-derivative} cannot directly be solved for $L$,
but an iterative method must be employed.
We treat the nonlinear constraint with an optimization on the Grassmann manifold in the next section and propose a variant of the conjugate gradient method for this purpose.

\subsection{Grassmann manifold}
\begin{figure}
  \centering
  \includegraphics[width=0.4\columnwidth]{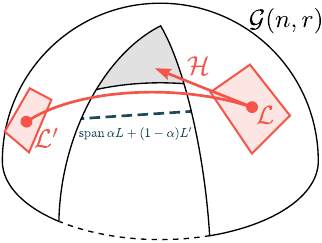}
  \caption{Illustration of the Grassmann manifold $\mathcal{G}(n,r)$. The red dots mark elements $\mathcal{L}$ and $\mathcal{L}^\prime$. Their linear interpolation is the dotted blue line which leaves $\mathcal{G}(n,r)$. The geodesic with direction $\mathcal{H}$ is given as the solid red line, see~\eqref{eq:geodesic}.
  }
  \label{fig:manifold-sketch}
\end{figure}

Let $F(L)$ denote the cost function in~\eqref{eq:omdc-cost-l}.
Observe that $F(L)$ is invariant under a right multiplication by an orthogonal matrix $R\in\mathcal{O}_{r\times r}$
\begin{equation*}
  F(L) = F( LR ).
\end{equation*}
Since the constraint~\eqref{eq:constraint-l} is also invariant under this multiplication, i.e.\ $R^\intercal L^\intercal L R= I$,  
the optimization problem~\eqref{eq:omdc-cost-l} subject to~\eqref{eq:constraint-l} can be solved for the subspace $\mathcal{L}$ 
\begin{equation}
  \mathcal{L} = \text{span}\,L = \text{span}\,LR,
\end{equation}
which is equal to the subspace spanned by $LR$ for any $R\in\mathcal{O}_{r\times r}$. 
Thus, the search space for the optimization is the Grassmann manifold $\mathcal{G}(n,r)$, which by definition is the set of all $r$-dimensional subspaces embedded in $\mathbb{R}^{n}$.
The subsequent operations are carried out in the tangent space to $\mathcal{L}\in\mathcal{G}(n,r)$, see Fig.~\ref{fig:manifold-sketch}.
We equip this tangent space with the Euclidian norm obtaining a Riemannian manifold.

The instrumental tool for optimization on a manifold is the geodesic.
The geodesic is defined as the shortest path between two manifold elements that stays on the manifold.
An illustration is shown in Fig.~\ref{fig:manifold-sketch} where the geodesic is compared to the linear combination of two bases.
The geodesic can be described by its origin $\mathcal{L}$ and a tangent vector $\mathcal{H}$, where the tangent vector $\mathcal{H}$ is represented by a matrix $H\in\mathbb{R}^{n\times r}$.
Then, the geodesic can be parametrized by the exponential map
\begin{equation*}
  \mathcal{L^\prime}(t) =\mathrm{Exp}_{\mathcal{L}} \, t\mathcal{H},
  \quad t \in [0,\,1],
\end{equation*}
where the representation in matrices reads
\begin{equation}
  \label{eq:geodesic}
  \begin{aligned}
    \hat{U}
    \Sigma V^\intercal &= H,
    \qquad \hat{U}\in\mathbb{R}^{n \times r},\, \Sigma,V\in\mathbb{R}^{r\times r}\\
   L^\prime (t) &= L V \cos (t\Sigma) V^\intercal +
   \hat{U}
   \sin (t\Sigma)V^\intercal
  \end{aligned}
\end{equation}
and $\hat{U}\Sigma V^\intercal$ denotes the thin singular value decomposition of $H$.
An optimization method must generate a sequence of tangential directions $H_k$ from which the next iterate $L_{k+1}$ can be computed with~\eqref{eq:geodesic}.

The remainder of this section translates the techniques for optimization in linear spaces to the Grassmann manifold.
The nonlinear conjugate gradient method~\cite{edelman1998geometry,Reineking2023a} proves to be a good choice, because it requires only first derivatives and still achieves superlinear convergence~\cite[sec.~5.2.2]{nocedal1999numerical}.
The matrix representation of gradient $G$ of the function $F$ with respect to the subspace $\mathcal{L}$ is defined as
\begin{equation}
  \label{eq:manifold-gradient}
  G = \frac{\partial F}{\partial\mathcal{L}} = \frac{\partial F}{\partial{L}} - L L^\intercal \frac{\partial F}{\partial{L}}
\end{equation}
where $\frac{\partial F}{\partial{L}}$ denotes the componentwise partial derivative of $F$ with respect to the orthogonal matrix $L$ from~\eqref{eq:cwise-derivative}.

The search direction $H_k$ of the $k$th step is combined from the steepest descend direction, i.e. the negative gradient $-G_k$, and the parallel-transported previous search direction $\tau H_{k-1}$ weighted with the conjugacy factor $\gamma_k$
\begin{equation}
  \label{eq:search-direction}
  H_k = -G_k + \gamma_k \, \tau H_{k-1}
\end{equation}
The search direction exists in the tangent space to $\mathcal{L}\in\mathcal{G}(n,r)$, as illustrated in Fig.~\ref{fig:manifold-sketch}.
Because the tangent spaces at two different points on the manifold differ, the parallel transport of the search direction $H_{k-1}$ along the geodesic to the next iterate on the manifold is necessary.
It reads~\cite{edelman1998geometry}
\begin{equation}
  \label{eq:parallel-transport}
  \begin{aligned}
    \tau H_{k-1} &= (-L_{k-1} V \sin(t\Sigma) +
    \hat{U}
    \cos(t\Sigma)) \Sigma V^\intercal\\
    \tau G_{k-1} &= G_{k-1} - (L_{k-1} V \sin t\Sigma\\
    &\quad\quad\quad\quad\quad\quad+ 
    \hat{U}
    (I - \cos t\Sigma))
    \hat{U}^\intercal G_{k-1}.
  \end{aligned}
\end{equation}
There are several choices for the conjugacy factor $\gamma_k$~\cite[sec.~5.2]{nocedal1999numerical}.
Numerical experiments have shown that the factor after Pollak-Ribiere for the Grassmann manifold~\cite[eq.~2.77]{edelman1998geometry}
\begin{equation}
  \gamma_k = \frac{\text{tr} \big( (G_k - \tau G_{k-1})^\intercal G_k \big)}
  {\text{tr} \big( G_k^\intercal G_k \big)}
\end{equation}
yields good convergence behavior.

The direction $H_k$ from~\eqref{eq:search-direction} defines a geodesic emanating from $L_k$ according to~\eqref{eq:geodesic}.
The step in direction $H_k$ to the next iterate $L_{k+1}$ is determined by minimizing the objective $F(L_k(t))$ in the scalar parameter $t$.
This can be done with standard line search methods.
The conjugate gradient method on the Grassmann manifold is summarized in Algorithm~\ref{alg:grassmann-cg}.
\begin{algorithm}[b]
  \KwResult{local minimizer $L^\ast$ for function $F$}
  Given $\epsilon, L_0$: $H_0 = -G_0$; $k \leftarrow 0$\;
  \While{$F(L_k) > \epsilon$}{
    solve line search $t^\ast = \arg \min F(L_k(t))$ along geodesic acc. to~\eqref{eq:geodesic}\;
    set $L_{k+1} \leftarrow L_k(t^\ast)$\;
    compute new derivative w.r.t. $L_{k+1}$ acc. to~\eqref{eq:cwise-derivative}\;
    obtain new manifold gradient $G_{k+1}$ from~\eqref{eq:manifold-gradient}\;
    transport $G_k$ and $H_k$ to new iterate with~\eqref{eq:parallel-transport}\;
    calculate new search direction $H_{k+1}$ with~\eqref{eq:search-direction}\;
    $k \leftarrow k + 1$\;
  }
  $L^\ast \leftarrow L_k$
  \caption{Conjugate gradient on Grassmann manifold~\cite[Alg.~13]{absil2008optimization}\label{alg:grassmann-cg}}
\end{algorithm}

\subsection{Implementation}
An implementation of the Grassmann conjugate gradient method and the derivative~\eqref{eq:cwise-derivative} results in a complexity of $\mathcal{O}(n)$ in every iteration of Algorithm~\ref{alg:grassmann-cg}, which may be prohibitive.
The mode matrix $L$ can, however, be expressed by a smaller matrix $\bar{L}\in\mathbb{R}^{q\times{r}}$ with the help of an orthonormal matrix $\bar{Q}\in\mathbb{R}^{n\times{q}}$
\begin{equation}
  \label{eq:mode-qr}
  L = \bar{Q}\bar{L}, \quad
  L^\intercal L = \bar{L}^\intercal\bar{Q}^\intercal\bar{Q}\bar{L}
  = \bar{L}^\intercal\bar{L} = I
\end{equation}
for some $q>r$.
We stress that the decomposition~\eqref{eq:mode-qr} and the steps to follow do not introduce an approximation.
In fact, it follows from~\eqref{eq:geodesic} that also the geodesic $L^\prime(t)$ can be expressed in terms of $\bar{Q}$ if the direction can be represented as $H=\bar{Q}\bar{H}$.
In the case of the conjugate gradient method, this is ensured if the gradient on the manifold is $G=\bar{Q}\bar{G}$.
The expression in terms of $\bar{Q}$ carries over to the derivative of the cost function $F(L)$ with respect to the modes $L$ via~\eqref{eq:manifold-gradient}.

It is clear from~\eqref{eq:cwise-derivative} that the derivative can be rewritten as
\begin{equation*}
  \frac{\partial F}{\partial L} = Y C_1 + X C_2
\end{equation*}
for some matrices $C_1,C_2\in\mathbb{R}^{(m-1)\times r}$.
Let the thin QR factorization~\cite[Th.~5.2.3]{golub2013matrix} of the snapshot matrix $S$ be denoted by
\begin{equation}
  S = \bar{Q}\bar{R}, \quad \bar{R} = [ \cdots r_k \cdots] \in \mathbb{R}^{m\times m}
\end{equation}
where $\bar{Q}\in\mathbb{R}^{n\times m}$, i.e. $q=m>r$ is fulfilled, and $\bar{Q}$ is orthonormal.
Then, the matrices $X$ and $Y$ can be put as
\begin{equation}
  \label{eq:xy-qr}
  \begin{aligned}
    X &= \bar{Q} \bar{X}, \quad &\bar{X} = [r_0 \cdots r_{m-2}] \\
    Y &= \bar{Q} \bar{Y}, \quad &\bar{Y} = [r_1 \cdots r_{m-1}]
  \end{aligned}
\end{equation}
which implies that the gradient of the cost function with respect to the modes reads
\begin{equation}
  \frac{\partial{F}}{\partial{L}} = \bar{Q} (\bar{C}_1 + \bar{C}_2)
\end{equation}
for some $\bar{C}_1,\bar{C}_2\in\mathbb{R}^{m\times r}$.
In summary, the conjugate gradient method for the given problem will always yield a mode matrix $L$ that can be represented in terms of the matrix $\bar{Q}$ generated from, e.g., the thin QR factorization of the snapshot matrix.
Thus, the OMDc problem can be solved equivalently on a $m$-dimensional space, by substituting~\eqref{eq:mode-qr} and~\eqref{eq:xy-qr} in~\eqref{eq:def-omdc-prob}.
The resulting problem reads
\begin{equation}
  \label{eq:def-omdc-prob-red}
  \min_{\bar{L},M,P} \Vert
  \bar{Y} - \bar{L}M\bar{L}^\intercal\bar{X} - \bar{L}P U\Vert_\text{F}.
\end{equation}
Note that the matrices $M$ and $P$ are identical to the original problem~\eqref{eq:def-omdc-prob}.
The complexity for the conjugate gradient method for the solution of~\eqref{eq:def-omdc-prob-red} is $\mathcal{O}(m)$ in every iteration.
This is a significant improvement whenever the number of snapshots $m$ is significantly smaller than the spatial resolution $n$.
We stress that the thin QR factorization has complexity $\mathcal{O}(n)$~\cite[Alg.~5.2.1]{golub2013matrix}, and, hence, is suited even for large problems because it is evaluated only once.
Alternatively, the snapshot matrix could be factorized with a thin singular value decomposition, also which has complexity $\mathcal{O}(n)$~\cite[Alg.~8.6.2]{golub2013matrix}, but still is significantly more costly than the thin QR factorization.
In the case of OMDc, we can assume $r < \text{rank}\, S \le m $ for the typical identification problem.
Thus, the rank-revealing thin singular value decomposition does not have any benefit here.
Nevertheless, the first $r$ left singular vectors of $\bar{R}$ can be computed with $\mathcal{O}(m)$ and can be used for initialization of $\bar{L}_0$ in Algorithm~\ref{alg:grassmann-cg}.

\section{Example}
\label{sec:woodchip}
\begin{figure}
  \centering
  \includegraphics{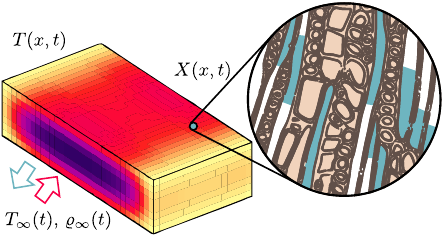}
  \caption{The wood chip is surrounded by hot dry air with temperature $T_\infty(t)$ and water vapor density $\varrho_\infty(t)$. A temperature distribution $T(x,t)$ and a moisture distribution $X(x,t)$ describe the interior the wood chip. During drying, the wood chip absorbs heat from the surrounding air and emits vapor. The pore structure, that actually exists in a wood particle, which is partially filled with water, is homogenized with spatially continuous quantities for the numerical treatment.}
  \label{fig:woodchip-sketch}
\end{figure}
We apply the method presented in section~\ref{sec:subspace-optimization} to a continuum model for the drying of wood chips.
Industrial drying processes are often carried out in rotating drums that may contain thousands or even tens of thousands of particles at a time.
In spite of this large number of particles, the processes in a single particle may still need to be resolved spatially because the particle is ``thermally thick'', i.e. it is too large to assume heat and water transport to be constant across its lengths
(for detailed simulations of industrial-scale processes see~\cite{scherer2016coupled}).
A computationally efficient model for the simulation of a single wood chip is obviously crucial in this setting~\cite{Reineking2023a}.
We use OMDc proposed in the present paper to derive such a model. 

The single particle drying process can be described with coupled diffusion equations for the 
the temperature 
$T(x,t)$ and  moisture $X(x,t)$ inside the wood chip particle, 
where $x$ refers to a point in space in this section and $t$ denotes time (see Fig.~\ref{fig:woodchip-sketch}).
The moisture $X(x,t)$ is defined as the ratio of the mass of the water contained in the wood chip and the dry mass of the wood chip.
Water contained in the wood chip evaporates to a surrounding stream of hot air. 
Although the model does not resolve the actual pore structure of wood, it accounts for (i) the anisotropic material structure due to the wood fibers, (ii) the dependencies of transport properties on the local temperature and moisture, and (iii) the evaporation of water which is confined to wood chip surface in this model.
Details on the continuum model are given in the appendix.

We assume the velocity of the surrounding air stream to be constant here, 
which implies constant heat and mass transfer coefficients.
The considered wood chip measures $5\,\text{mm} \times 10\,\text{mm} \times 20\,\text{mm}$ in height, width and length, respectively.
The partial differential equations for $X(x,t)$ and $T(x,t)$ are spatially discretized with 8000 finite-volume cells.
The finite-volume simulation is carried out for $1250\,\text{s}$ with a time step of $0.1\,\text{s}$, beginning with an initial moisture and temperature of $X(x,0) = 0.8$ and $T(x,0) = 298.15\,\mathrm{K}$, respectively, i.e., a moist wood chip at nearly room temperature is dried.
The simulation results for $X$ and $T$ are collected every $\Delta t = 12.5\,\text{s}$ into the snapshot matrix
\begin{equation}
  \label{eq:woodchip-snaps}
  S = \begin{bmatrix}
    \cdots & X(x,t_k) & \cdots \\
    \cdots & T(x,t_k) & \cdots \\
  \end{bmatrix}, \quad k = 0,\ldots,100
\end{equation}
which implies $n=16000$.
Including the initial state, $m=101$ snapshots are collected.

Two inputs are available to control the process: the temperature $T_\infty(t)$ and water vapor density $\varrho_{\infty}(t)$ of the air stream.
They are chosen to be 
\begin{equation}
  \label{eq:def-ambient-inputs}
  \begin{aligned}
    T_\infty(t) &= 375\,\text{K}, \, 0 \le t \le 1250 \\
    \varrho_\infty(t) &=
    \begin{cases}
      0.0350\,\text{kg}/\text{m}^3, & \quad0 \le t < 100 \\
      0.0175\,\text{kg}/\text{m}^3, & 100 \le t < 200 \\
      0.0070\,\text{kg}/\text{m}^3, & 200 \le t \le 1250
    \end{cases}
  \end{aligned}
\end{equation}
in the example treated here.
These inputs are collected in the input matrix $U$
\begin{equation}
  \label{eq:woodchip-inputs}
  U = \begin{bmatrix}
    \cdots T_\infty(t_k) \cdots \\
    \cdots \varrho_\infty(t_k) \cdots
  \end{bmatrix}, \quad k=0,\ldots,99
\end{equation}
Observe that~\eqref{eq:def-ambient-inputs} results in two independent rows of $U$. 

\begin{figure}
  \centering
  \includegraphics[trim=0pt 4pt 0pt 2pt,clip]{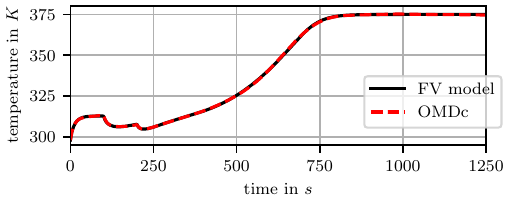}
  \caption{Mean temperature of the wood chip over time.}
  \label{fig:temperature-response}
\end{figure}

\begin{figure}
  \centering
  \includegraphics[trim=0pt 4pt 0pt 2pt,clip]{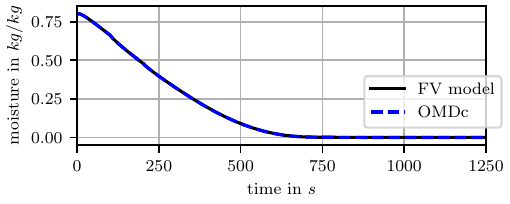}
  \caption{Mean moisture of the wood chip over time.}
  \label{fig:moisture-response}
\end{figure}

The results of the finite-volume simulation are shown in Figs.~\ref{fig:temperature-response} and~\ref{fig:moisture-response}.
In the initial $100\,\text{s}$, the average temperature of the wood chip increases to a plateau of approximately $312\,\text{K}$.
There, the heat transfer from the surrounding air to the wood chip is in equilibrium with the heat consumption due to the evaporation of moisture.
Then, the ambient water vapor density $\varrho_\infty(t)$ is decreased at $t=100\,\text{s}$.
As a consequence, the evaporation rate increases, resulting in a greater heat consumption and a reduction of the average wood chip temperature. 
The same effect can be observed at $t=200\,\mathrm{s}$, where the second decrease of the ambient water vapor density takes place.
Afterwards, the wood chip heats up to the ambient temperature of $375\,\text{K}$.

\begin{figure}
  \centering
  \includegraphics[trim=0pt 4pt 0pt 2pt,clip]{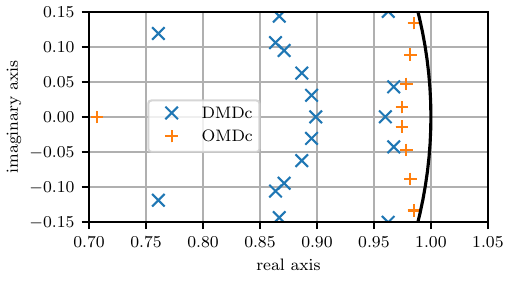}
  \caption{Eigenvalues of the identified linear system for DMDc and OMDc in the complex plane. The unit circle is drawn in black for reference.}
  \label{fig:eigenvalues}
\end{figure}

The snapshots~\eqref{eq:woodchip-snaps} and inputs~\eqref{eq:woodchip-inputs} are processed as described in section~\ref{sec:subspace-optimization}.
It is known from a previous study~\cite{Reineking2023a} that $5$ modes for each temperature and moisture are sufficient for a POD-Galerkin reduced model.
The dimension for OMDc was therefore chosen to be $r=10$ for comparison.
A good agreement of the linear reduced model from OMDc and the finite-volume simulation is evident from Figs.~\ref{fig:temperature-response} and~\ref{fig:moisture-response}.
Furthermore, the eigenvalues of the system matrices computed with OMDc and DMDc differ, as shown in Fig.~\ref{fig:eigenvalues}.
This implies that the computed modes $L^\ast$ are not a subset of the DMDc modes $\hat{\Phi}_r$.

A comparison of the computation times of the presented methods is given in Tab.~\ref{tab:computation-times}.
All computations were carried out on the same computer with an Intel i7-10700K CPU running at 3.8GHz using only a single thread.
Note that the finite-volume simulation uses a different time step than the reduced model, resulting in a by far larger total computation time. For the given example, Algorithm~\ref{alg:grassmann-cg} converged after $600$ iterations. Still, it is only one order of magnitude slower than DMDc, which is not an iterative method.
We stress that this relative difference decreases for a larger ratio of $n$ and $m$, as the optimization problem~\eqref{eq:def-omdc-prob-red} is independent of $n$.
\begin{table}
  \centering
  \caption{summary of computation times}
  \label{tab:computation-times}
  \begin{tabular}{lrr}
    \hline
    Computation & per iteration & total \\
    \hline
    finite-volume simulation & $4.75\,\text{ms}$    & $59410.2\,\text{ms}$ \\
    DMDc    & $-$ & $487.2\,\text{ms}$ \\
    OMDc & $8.73\,\text{ms}$ & $5509.1\,\text{ms}$\\
    OMDc reduced model & $0.014\,\text{ms}$ & $1.4\,\text{ms}$\\
    \hline
  \end{tabular}
\end{table}

\section{Conclusion}\label{sec:conclusions}
We presented a new variant of dynamic mode decomposition, optimal mode decomposition for control, which extends optimal mode decomposition from autonomous systems to systems with inputs.
The application of the thin QR factorization renders the implementation competitive 
even though an optimization on a manifold must be carried out.

We illustrated the proposed method by applying it to coupled PDEs that model the drying of wood chips. The resulting reduced model is ideally suited for use in simulations of industrial drying systems, where PDE models are too computationally expensive because thousands or tens of thousands of particles must be modeled~\cite{scherer2016coupled}. 

In our future work, we will investigate parametrized snapshot data.
It has been shown in~\cite{Mjalled2023b} that parametrized neural networks for the state prediction combined with an interpolation of POD modes on the Grassmann manifold yields good results.
We will treat the interpolation of the identified linear systems in reduced space and couple it to the interpolation of modes on the Grassmann manifold to create parametrized reduced models.

\appendix
\section{PDE Model}
\label{app:pde-model}
The temperature distribution $T(x,t)$ and moisture distribution $X(x,t)$ inside a wood chip particle are described by coupled partial differential equations
\begin{equation}
  \label{eq:pde-drying}
  \begin{aligned}
    \rho(X) c_p(X) \dot{T}(x,t) &= \nabla \cdot (\lambda(X) \cdot \nabla T(x,t)) \\
    \dot{X}(x,t) &= \nabla \cdot (\delta(T)\cdot \nabla X(x,t)).
  \end{aligned}
\end{equation}
The heat and moisture transport are modelled by conduction and diffusion, respectively.
The model is nonlinear as the conductivity $\lambda(X)$ and the diffusivity $\delta(T)$ depend on the local temperature and moisture.
Further, conductivity and diffusivity are tensors accounting for anisotropic material structure of wood fibers.
The model is closed by
\begin{equation}
  \begin{aligned}
    \delta(X) \cdot X \cdot \vec{n} &= \beta \left(\varrho_\infty - \varrho\right)
    =\dot{m}_\mathrm{w} \\
    \lambda(X) \cdot \nabla T \cdot \vec{n} &= \alpha \, (T_\infty -  T) + \Delta h(X,T) \,\beta \, \dot{m}_\mathrm{w}
    .
  \end{aligned}
\end{equation}
This model confines the evaporation of water to the wood chip surface, leading to a nonlinear boundary term.
The complete model and its verification can be found in~\cite{scherer2016coupled}.

\begin{table}
  \centering
  \caption{Properties of the dry wood chip according to~\cite{Reineking2023a}}
  \begin{tabular}{lcrl}
    \hline
    Property & symbol & value~~~~&  \\
    \hline
    density       & $\rho_\text{dry}$    & $500$&$\text{kg}\,\text{m}^{-3}$ \\
    heat capacity    & $c_{p,\text{dry}}$   & $1500$&$\text{J}\,\text{kg}^{-1}\text{K}^{-1}$ \\
    thermal conductivity  & $\lambda_\text{dry}$ & $0.12$&$\text{W}\,\text{m}^{-1}\text{K}^{-1}$ \\
    effective diffusivity           & $\delta_\text{eff}$           & $2\cdot 10^{-9}$&$\text{m}^2\text{s}^{-1}$ \\
    heat transfer coeff.  & $\alpha$             & $45$&$\text{W}\,\text{m}^{-2}\text{K}^{-1}$ \\
    mass transfer coeff.  & $\beta$              & $0.075$&$\text{m}\,\text{s}^{-1}$\\
    \hline
  \end{tabular}
\end{table}

\bibliographystyle{IEEEtran}
\bibliography{references}

\end{document}